\documentclass{article}
\usepackage{amssymb}
\usepackage{amsfonts}
\usepackage{amsmath}

\input{tcilatex}
\begin{document}

\title{Approximation by Choquet Integral Operators}
\author{Sorin G. Gal \\
Department of Mathematics and Computer Science, \\
University of Oradea, \\
Universitatii 1, 410087, Oradea, Romania\\
E-mail: \textit{galso@uoradea.ro}}
\date{}
\maketitle

\begin{abstract}
The main aim of this paper is to show that the nonlinear Choquet integral can be used to construct nonlinear approximation operators, exactly as by the use in probability of the Lebesgue-type integral, linear and positive approximation operators are constructed.
The so-called Feller constructive scheme is generalized, by introducing discrete and non-discrete nonlinear approximation operators in terms of the nonlinear Choquet integral with respect to a monotone and subadditive set function.
As particular cases, Bernstein-Choquet and Picard-Choquet operators are introduced, for which qualitative and quantitative approximation properties are obtained. In some subclasses of functions, they have better approximation properties than the classical Berstein and Picard operators.
\end{abstract}



\section{Introduction}

A well-known general scheme in constructing linear and positive approximation operators, is the Feller's probabilistic scheme
(see \cite{Feller}, Chapter 7, or e.g. \cite{Alto}, Section 5.2, pp. 283-319), by which
to any continuous and bounded function $f:\mathbb{R}\to \mathbb{R}$, approximation operators of the form
\begin{equation}\label{int1}
L_{n}(f)(x)=\int_{\Omega}f\circ Z(n, x) d P=\int_{\mathbb{R}}f d P_{Z(n, x)},
\end{equation}
are attached, where $P$ is a probability on the measurable space $(\Omega, {\cal{C}})$, $Z:\mathbb{N}\times I \to {\cal{M}}_{2}(\Omega)$, with
$I$ a subinterval of $\mathbb{R}$, ${\cal{M}}_{2}(\Omega)$ represents the space of all random variables whose square is integrable
on $\Omega$ with respect to the probability $P$ and $P_{Z(n, x)}$ denotes the distribution of the random variable $Z(n, x)$ with respect to $P$ defined by $P_{Z(n, x)}(B)=P(Z^{-1}(n, x))$, for all Borel measurable subsets of $\mathbb{R}$. Then, denoting by $E(Z(n, x))$ and $Var(Z(n, x))$ the expectance and the variance of the random variable $Z(n, x)$, respectively, and supposing that $\lim_{n\to \infty}E(Z(n, x))=x$,
$\lim_{n\to \infty}Var(Z(n, x))=0$, uniformly on $I$, it is proved that for all $f$ as above, $L_{n}(f)$ converges to $f$ uniformly on each compact subinterval of $I$.

In addition, if for the random variable $Z(n, x)$ its probability density function $\lambda_{n, x}$ is known, then for any $f$ we can write
\begin{equation}\label{int2}
L_{n}(f)(x)=\int_{\mathbb{R}}f d P_{Z(n, x)}=\int_{\mathbb{R}} f(t)\cdot \lambda_{n, x}(t) d P(t),
\end{equation}
formula which is useful for the concrete construction and quantitative estimations in the approximation by the operators $L_{n}(f)(x)$.

It is worth noting that for the most classical linear and positive operators, the approximation properties can be studied by using
their representations in the above forms (\ref{int1}) and (\ref{int2}) (see e.g. \cite{Alto}, Section 5.2, pp. 283-319).

The main aim of this paper is to generalize the constructive formulas in (\ref{int1}) and (\ref{int2}) for the case when $P$ is a monotone, subadditive set function, not necessarily (countable) additive as in the case of the probability.

More precise, we consider a Feller kind scheme based on the Choquet integral with respect to a monotone and subadditive set function (sometimes called capacity), for the construction of approximation sequences of nonlinear operators. Since in the case when the set function $P$ is monotone and subadditive, in general the formula (\ref{int2}) does not hold, the approximation sequences of nonlinear operators can formally be constructed by two distinct ways :  $L_{n}(f)(x)$ given by the formula in (\ref{int1}) and $L_{n}(f)(x)$ given by the last integral on the right-hand side of (\ref{int2}).
However, it is worth noting that according to the Radon-Nikodym type result for the Choquet integral in \cite{Graf} (see also e.g. \cite{Mol}, p. 75, Theorem 5.10),
for very special subclasses of monotone and subadditive set functions, the formulas for $L_{n}(f)(x)$ in (\ref{int2}) still holds.

Bernstein-Choquet and Picard-Choquet type operators are introduced and their qualitative and quantitative approximation properties are studied. In some subclasses of functions, they have better approximation properties than the classical Bernstein and Picard operators.

Note that similar approaches were considered in \cite{Gal-Pos} and \cite{Gal-Pos2} for the study of the nonlinear approximation operators constructed in terms of the possibilistic integral, including the important classes of the so-called max-product approximation operators.

\section{Preliminaries}

The main aim of this section is to present known concepts and results used in the next section.

{\bf Definition 2.1.} Let $(\Omega, {\cal{C}})$ be a measurable space, i.e. $\Omega$ is a nonempty set and ${\cal{C}}$ be a $\sigma$-ring (or $\sigma$-algebra) of subsets in $\Omega$ with $\emptyset\in {\cal{C}}$.

(i) (see e.g. \cite{WK1}, p. 63) The set function $\mu:{\cal{C}}\to [0, +\infty]$ is called a monotone set function (or capacity) if
$\mu(\emptyset)=0$ and $A, B\in {\cal{C}}$, with $A\subset B$, implies $\mu(A)\le \mu(B)$. If $\Omega\in {\cal{C}}$, then
$\mu$ is called normalized if $\mu(\Omega)=1$.

(ii) (see \cite{Choquet}, or e.g. \cite{WK1}, p. 179) Let $\mu$ be a normalized, monotone set function and denote by ${\cal{G}}$ the class of all finite nonnegative functions defined and measurable on $(\Omega, {\cal{C}})$. Recall that $X:\Omega\to \mathbb{R}$ is measurable (or more precisely ${\cal{C}}$-measurable) if for any $B$, Borelian subset in $\mathbb{R}$, we have $X^{-1}(B)\in {\cal{C}}$.

For $A\in {\cal{C}}$ and $X\in {\cal{G}}$, the Choquet integral of $X$ on $A$ with respect to a monotone set function $\mu$ is defined by
$$(C)\int_{A}X d\mu=\int_{0}^{\infty}\mu(F_{\alpha}(X)\bigcap A)d\alpha,$$
where $F_{\alpha}(X)=\{\omega\in \Omega; X(\omega)\ge \alpha\}$. If $(C)\int_{A}X d\mu <+\infty$ then $X$ is called Choquet integrable on $A$.

If $X:\Omega \to \mathbb{R}$ is of arbitrary sign, then the Choquet integral is defined by (see \cite{WK1}, p. 233)
$$(C)\int_{A} X d\mu=\int_{0}^{+\infty}\mu(F_{\alpha}(X)\bigcap A)d\alpha+\int_{-\infty}^{0}[\mu(F_{\alpha}(X)\bigcap A)-\mu(A)]d \alpha.$$

When $\mu$ is the Lebesgue measure (i.e. countably additive), then the Choquet integral $(C)\int_{A}X d\mu$ reduces to the Lebesgue integral.

In the case when $\Omega=\{\omega_{0}, \omega_{1}, ..., \omega_{n}\}$ is a finite set and $\mu:{\cal{P}}(\Omega)\to \mathbb{R}_{+}$ is a monotone set function, then the Choquet integral of $X:\Omega\to \mathbb{R}$ with respect to $\mu$ reduces to the formula (see e.g. \cite{WK1}, p. 239)
$$(C)\int_{\Omega}X d\mu=\sum_{i=0}^{n}[X(\omega_{i}^{*})-X(\omega_{i-1}^{*})]\mu(\{\omega_{i}^{*}, \omega_{i+1}^{*}, ..., \omega_{n}^{*}\})$$
$$=\sum_{i=0}^{n}X(\omega_{i}^{*})[\mu(\{\omega_{i}^{*}, ..., \omega_{n}^{*}\})-\mu(\{\omega_{i+1}^{*}, ..., \omega_{n}^{*}\})],$$
where $(\omega_{0}^{*}, \omega_{1}^{*}, \omega_{2}^{*}, ..., \omega_{n}^{*})$ is a permutation of $(\omega_{0}, \omega_{1}, \omega_{2}, ..., \omega_{n})$ such that
$$X(\omega_{0}^{*})\le X(\omega_{1}^{*})\le X(\omega_{2}^{*})\le ... \le X(\omega_{n}^{*})$$
and by convention, in the first sum we take $X(\omega_{-1}^{*})=0$, while in the second sum we take $\mu(\{\omega_{n+1}^{*}, \omega_{n}^{*}\})=0$.

It is worth mentioning that the above concept of Choquet integral appears for the first time in the 1925 paper of Vitali \cite{Vitali} (see \cite{WK1}, p. 245).

(iii) Given a measurable space $(\Omega, {\cal{C}})$ and $\mu:{\cal{C}}\to [0, +\infty)$ a normalized, monotone set function, for each ${\cal{C}}$-measurable
$X:\Omega\to \mathbb{R}$, we can define its distribution (as a set function) with respect to $\mu$, by the formula
$$\mu_{X}:{\cal{B}}\to \mathbb{R}_{+}, \, \mu_{X}(B)=\mu(X^{-1}(B))=\mu(\{\omega\in \Omega ; X(\omega)\in B\}),\, B\in {\cal{B}},$$
where ${\cal{B}}$ is the class of all Borel measurable subsets in $\mathbb{R}$.
It is clear that $\mu_{X}$ is a normalized, monotone set function on ${\cal{B}}$.

We say that $X$ (and its distribution $\mu_{X}$) have the Choquet density $\lambda_{X}:\mathbb{R}\to \mathbb{R}_{+}$ (with respect to $\mu$), if $\lambda_{X}$ is ${\cal{B}}$-measurable, $(C)\int_{\mathbb{R}}\lambda_{X} d \mu=1$ and
$$\mu_{X}(B)=\mu(X^{-1}(B))=(C)\int_{B}\lambda_{X} d\mu=(C)\int_{B} d \mu_{X},\, \mbox{ for all } B\in {\cal{B}}.$$

(iv) Given $X:\Omega\to \mathbb{R}$ and $f:\mathbb{R}\to \mathbb{R}$, ${\cal{B}}$-measurable, the Choquet integral of $f$ with respect to the distribution $\mu_{X}$ is denoted by $(C)\int_{\mathbb{R}}f(t) d\mu_{X}(t)$.

(v) For a measurable $X:\Omega \to \mathbb{R}$, its Choquet expectance is defined by
$E_{Ch}(X)=(C)\int_{\Omega}X d\mu$. Also, its Choquet variance is defined by $VAR_{Ch}(X)=E_{Ch}((X-E_{Ch}(X))^{2})$.

The following result is important for our considerations and seems to be new.

{\bf Lemma 2.2.} {\it If $\mu:{\cal{C}}\to \mathbb{R}_{+}$ is a normalized, monotone set function, $X$ is ${\cal{C}}$-measurable and $f:\mathbb{R}\to\mathbb{R}$ is ${\cal{B}}$-measurable, then we have
$$(C)\int_{\mathbb{R}}f(t) d\mu_{X}(t)=(C)\int_{\Omega}f \circ X d\mu$$
$$=\int_{0}^{+\infty}\mu(F_{\alpha}(f \circ X))d \alpha+\int_{-\infty}^{0}[\mu(F_{\alpha}(f \circ X))-1]d \alpha.$$}
{\bf Proof.} Denoting $B_{\alpha}=\{t\in \mathbb{R} ; f(t)\ge \alpha\}$, it follows
$$X^{-1}(B_{\alpha})=\{\omega\in \Omega ; X(\omega)\in B_{\alpha}\}=\{\omega\in \Omega ; f[X(\omega)]\ge \alpha\},$$
which immediately implies
$$(C)\int_{\mathbb{R}}f(t) d\mu_{X}(t)=\int_{0}^{+\infty}\mu_{X}(B_{\alpha})d \alpha + \int_{-\infty}^{0}[\mu_{X}(B_{\alpha})-1]d \alpha$$
$$=\int_{0}^{+\infty}\mu(X^{-1}(B_{\alpha}))d \alpha + \int_{-\infty}^{0}[\mu(X^{-1}(B_{\alpha}))-1] d \alpha
=(C)\int_{\Omega}f \circ X d\mu.$$
Note that more general, for any $T\in {\cal{B}}$, denoting $\Omega_{T}=X^{-1}(T)$, we have
$$(C)\int_{T}f(t) d\mu_{X}(t)=(C)\int_{\Omega_{T}}f \circ X d\mu$$
$$=\int_{0}^{+\infty}\mu(F_{\alpha}(f \circ X)\bigcap \Omega_{T})d \alpha + \int_{-\infty}^{0}[\mu(F_{\alpha}(f \circ X)\bigcap \Omega_{T})-\mu(\Omega_{T})]d \alpha$$
and the lemma is proved. $\hfill \square$

In what follows, we list some know properties we need for our reasonings.

{\bf Remark 2.3.} Let us suppose that $\mu$ is a monotone set function (sometimes such a $\mu$ is called capacity). Then, the following properties hold :

(i) $(C)\int_{A}$ is non-additive (i.e. $(C)\int_{A}(f+g)d\mu \not= (C)\int_{A}f d\mu + (C)\int_{A}g d \mu$) but it is positive homogeneous, i.e.
for all $a\ge 0$ we have $(C)\int_{A}af d\mu = a\cdot (C)\int_{A}f d\mu$ (for $f\ge 0$ see, e.g., \cite{WK1}, Theorem 11.2, (5), p. 228 and for $f$ of arbitrary sign, see e.g., \cite{Denn}, p. 64, Proposition 5.1, (ii)).

If $f\le g$ on $A$ then the Choquet integral is monotone, that is
$(C)\int_{A}f d\mu \le (C)\int_{A}g d\mu$ (see e.g. \cite{WK1}, p. 228, Theorem 11.2, (3) for $f, g\ge 0$ and p. 232 for $f, g$ of arbitrary sign).

If $\mu$ is submodular too (i.e. $\mu(A\bigcup B)+\mu(A\bigcap B)\le \mu(A)+\mu(B)$ for all $A, B$) then the Choquet integral is subadditive, that is $(C)\int_{A}(f + g) d\mu \le (C)\int_{A}f d\mu + (C)\int_{A}g d\mu$, for all $f, g$ of arbitrary sign (see e.g. \cite{Denn}, p. 75, Theorem 6.3).

If $\overline{\mu}$ denotes the dual set function of $\mu$ (that is $\overline{\mu}(A)=\mu(\Omega)-\mu(\Omega\setminus A)$, for all $A\in {\cal{C}}$), then for all $f$ of arbitrary sign we have $(C)\int_{A}(-f) d\mu=-(C)\int_{A} f d \, \overline{\mu}$ (see e.g. \cite{WK1}, Theorem 11.7, p. 233).

If $c\in \mathbb{R}$ and $f$ is of arbitrary sign, then $(C)\int_{A}(f+c)d \mu = (C)\int_{A}f d\mu + c\cdot \mu(A)$ (see e.g. \cite{WK1}, pp. 232-233, or \cite{Denn}, p. 65).

By the definition of the Choquet integral, if $F\ge 0$ and $\mu$ is subadditive, then it is immediate that
$$(C)\int_{A\bigcup B}F d\mu \le (C)\int_{A}F d\mu + (C)\int_{B}F d\mu.$$
Note that if $\mu$ is submodular then it is clear that it is subadditive too.

(ii) Simple concrete examples of monotone and submodular set functions $\mu$, can be obtained from a probability measure $M$ on ${\cal{P}}(\mathbb{X})$ (i.e. $M(\emptyset)=0$, $M(\mathbb{X})=1$ and $M$ is countable additive), by the formula $\mu(A)=\gamma(M(A))$, where
$\gamma :[0, 1]\to [0, 1]$ is an increasing and concave function, with $\gamma(0)=0$, $\gamma(1)=1$ (see e.g. \cite{Denn}, pp. 16-17, Example 2.1).

Also, any possibility measure $\mu$ is monotone and submodular. While the monotonicity is immediate from the axiom  $\mu(A\bigcup B)=\max\{\mu(A), \mu(B)\}$, the submodularity is immediate from the property $\mu(A\bigcap B)\le \min\{\mu(A), \mu(B)\}$.

(iii) Many other properties of the Choquet integral can be found in e.g. Chapter 11 in \cite{WK1}, or in \cite{Denn}.

The following Chebyshev-type inequality for the Choquet integral obtained directly from Corollary 3.1 in \cite{Wang} (see also \cite{She}) will also be useful.

{\bf Theorem 2.4.} (Chebyshev's inequality) {\it If $\Omega\in {\cal{C}}$, $\mu$ is a monotone set function and $F:\Omega\to \mathbb{R}$ is ${\cal{C}}$-measurable, then for any $r>0$ we have
$$\mu\left (\left \{s\in \Omega; \left |F(s)-(C)\int_{\Omega}F d\mu\right |\ge r\right \}\right )\le \frac{(C)\int_{\Omega}\left (F - (C)\int_{\Omega}F d\mu\right )^{2}d \mu}{r^{2}}.$$}

\section{Nonlinear Choquet Integral Operators}

The aim of this section is to introduce operators constructed in terms of the Choquet integral and to study their approximation properties.
By analogy to the Feller's random scheme in probability theory which produce linear and positive approximation operators, we will consider a similar approximation scheme, but which will produce nonlinear approximation operators in terms of the Choquet integral.

Given a measurable space $(\Omega, {\cal{C}})$, $\mu:{\cal{C}}\to [0, +\infty)$ a normalized, monotone set function and $I\subset \mathbb{R}$ a real interval (bounded or unbounded), let us consider a mapping $Z:\mathbb{N}\times I \to Mes_{{\cal{C}}}(\Omega)$, where $Mes_{{\cal{C}}}(\Omega)$ denotes the class of all $X:\Omega \to \mathbb{R}$ which are ${\cal{C}}$-measurable.

For any $(n, x)\in \mathbb{N}\times I$, denote
$$E_{Ch}(Z(n, x))=(C)\int_{\Omega} Z(n, x)d \mu:=\alpha_{n, x} \mbox{ and } VAR_{Ch}(Z(n, x)):=\sigma^{2}_{n, x}.$$
Now, according to the idea in the Feller's scheme, to continuous $f:\mathbb{R}\to \mathbb{R}$, let us attach a sequence of operators by the formula in Lemma 2.2,
$$L_{n}(f)(x):=(C)\int_{\mathbb{R}}f(t) d\mu_{Z(n, x)}(t)=(C)\int_{\Omega}f \circ Z(n, x) d\mu, \, x\in I, \, n\in \mathbb{N}.$$

{\bf Theorem 3.1.} {\it Let $\mu$ be a normalized, monotone and subadditive set function and $f:\mathbb{R}\to \mathbb{R}$ be uniformly continuous and bounded on $\mathbb{R}$.
Suppose that $\lim_{n\to +\infty}\alpha_{n, x}=x$, uniformly on $I$ and $\lim_{n\to +\infty}\sigma^{2}_{n, x}=0$, uniformly on each compact subinterval of $I$. Then, $\lim_{n\to \infty}L_{n}(f)=f$, uniformly on any compact subinterval of $I$.}

{\bf Proof.} First of all, it is worth noting that since $f$ is bounded, both integrals defining $L_{n}(f)(x)$ are finite, for all $x\in I$ and $n\in \mathbb{N}$. Indeed, denoting $H(\alpha)=\mu(\{\omega\in \Omega ; f[Z(n, x)(\omega)]\ge \alpha\})$, it is clear that $H(\alpha)$ is a finite (bounded) and nonincreasing function on $(-\infty, +\infty)$ and since $f$ is bounded, there exists $m<0$ and $M>0$, such that we can write
$$(C)\int_{\Omega}f\circ Z(n, x) d\mu=\int_{0}^{M}H(\alpha) d \alpha +\int_{m}^{0}[H(\alpha) - 1] d \alpha,$$
which obviously that is a finite number.

Now, since $f$ is uniformly continuous on $\mathbb{R}$, for any $\varepsilon >0$, there exist $\delta > 0$, such that for all $t, x\in \mathbb{R}$ with $|t-x| < \delta$, we have $|f(t)-f(s)|\le \varepsilon/2$.

Let $G:\mathbb{R}\to \mathbb{R}$ be bounded on $\mathbb{R}$ and ${\cal{B}}$-measurable. Since $\mu$ is subadditive
on ${\cal{C}}$, we get
$$\overline{\mu}(A)=\mu(\Omega)-\mu(\Omega \setminus A)\le \mu(A).$$
Now, applying the operator $L_{n}$ to the obvious inequality
$-|G|(x)\le G(x)\le |G|(x)$, for all $x\in \mathbb{R}$, we obtain
$$L_{n}(-|G|)(x)\le L_{n}(G)(x)\le L_{n}(|G|)(x), x\in \mathbb{R}.$$
But by the inequality $\overline{\mu}(A)\le \mu(A)$, for all $A\in {\cal{C}}$ in Remark 2.3, (i), and by Lemma 2.2, we get
$$L_{n}(-|G|)(x)=(C)\int_{\Omega}-|G|\circ Z(n, x) d \mu = - (C)\int_{\Omega}|G|\circ Z(n, x) d\, \overline{\mu}$$
$$=-\int_{0}^{+\infty}\overline{\mu}[F_{\alpha}(|G|\circ Z(n, x))]d \alpha - \int_{-\infty}^{0}(\overline{\mu}[F_{\alpha}(|G|\circ Z(n, x))]- 1)d \alpha$$
$$\ge - \int_{0}^{+\infty}\mu[F_{\alpha}(|G|\circ Z(n, x))]d \alpha - \int_{-\infty}^{0}(\mu[F_{\alpha}(|G|\circ Z(n, x))]- 1)d \alpha$$
$$=-(C)\int_{\Omega}|G|\circ Z(n, x) d\mu = - L_{n}(|G|)(x),$$
which leads to
$-L_{n}(|G|)(x)\le L_{n}(G)(x)\le L_{n}(|G|)(x)$, equivalent to
$$|L_{n}(G)(x)|\le L_{n}(|G|)(x), \mbox{ for all } x\in \mathbb{R}.$$
Above we used the relationship $\overline{\mu}(\Omega)=\mu(\Omega)=1$.

By using Theorem 11.6, p. 232 in \cite{WK1} (see also the last but one property in Remark 2.3, (i)) and the above property for $G(t):=f(t) - f(\alpha_{n, x})$, it follows
$$|L_{n}(f)(x)-f(\alpha_{n,x})|=\left |(C)\int_{\mathbb{R}}(f(t) - f(\alpha_{n, x}))d\mu_{Z(n, x)}(t)\right |$$
$$\le (C)\int_{\mathbb{R}}|f(t) - f(\alpha_{n, x})|d\mu_{Z(n, x)}(t).$$

Let us consider the decomposition
$$\mathbb{R}=\{t\in \mathbb{R} ; |t-\alpha_{n, x}|<\delta\} \bigcup \{t\in \mathbb{R} ; |t-\alpha_{n, x}|\ge \delta\}:=T_{1}\bigcup T_{2}$$
and denote $\Omega_{T_{1}}=Z^{-1}(n, x)(T_{1})$, $\Omega_{T_{2}}=Z^{-1}(n, x)(T_{2})$.

Since $\mu$ is subadditive, this easily implies the subadditivity of $\mu_{Z(n, x)}$, which implies
$$(C)\int_{\mathbb{R}}|f(t) - f(\alpha_{n, x})|d\mu_{Z(n, x)}(t)$$
$$=\int_{0}^{\infty}\mu_{Z(n, x)}(\{t\in T_{1}\bigcup T_{2} ; |f(t) - f(\alpha_{n, x})|\ge \alpha\})d \alpha$$
$$\le \int_{0}^{\infty}\left [\mu_{Z(n, x)}(\{t\in T_{1} ; |f(t) - f(\alpha_{n, x})|\ge \alpha\}) \right .$$
$$\left .+\mu_{Z(n, x)}(\{t\in T_{2} ; |f(t) - f(\alpha_{n, x})|\ge \alpha\})\right ]d \alpha$$
$$= \int_{0}^{\infty}\mu_{Z(n, x)}(\{t\in T_{1} ; |f(t) - f(\alpha_{n, x})|\ge \alpha\})d \alpha$$
$$+ \int_{0}^{\infty}\mu_{Z(n, x)}(\{t\in T_{2} ; |f(t) - f(\alpha_{n, x})|\ge \alpha\})d \alpha$$
$$=(C)\int_{T_{1}}|f(t) - f(\alpha_{n, x})|d\mu_{Z(n, x)}(t)
+(C)\int_{T_{2}}|f(t) - f(\alpha_{n, x})|d\mu_{Z(n, x)}(t)$$
$$\le \frac{\varepsilon}{2} (C) \int_{T_{1}} 1 \cdot d \mu_{Z(n, x)}(t)
+2 \|f\|\cdot \mu(\{|Z(n, x)-\alpha_{n, x}|\ge \delta\})$$
$$\le \frac{\varepsilon}{2}+2\|f\|\cdot \sigma^{2}_{n, x}\cdot \delta^{-2}
\le \frac{\varepsilon}{2}+\frac{\varepsilon}{2}=\varepsilon,$$
for all $n\ge n_{0}$, uniformly on any compact subinterval of $I$.

Above $\|f\|=\sup\{|f(t)| ; t\in \mathbb{R}\}$ and we also used the relationships (based on Remark 2.3, (i), first property)
$$(C) \int_{T_{1}} 1 \cdot d \mu_{Z(n, x)}(t)=(C)\int_{\Omega_{T_{1}}} 1\cdot d\mu\le (C)\int_{\Omega} 1\cdot d\mu=\mu(\Omega)=1,$$
$$(C)\int_{T_{2}} |f(t) - f(\alpha_{n, x})|d\mu_{Z(n, x)}(t)\le 2\|f\|\cdot (C)\int_{T_{2}} 1\cdot  d \mu_{Z(n, x)}(t)$$
$$=2\|f\|\cdot (C)\int_{\Omega_{T_{2}}}1 d \mu
=2\|f\|\cdot \mu(\Omega_{T_{2}})=2\|f\|\cdot \mu(\{\omega\in \Omega ; |Z(n, x)(\omega)-\alpha_{n, x}|\ge \delta\})$$
and the Chebyshev's inequality in Theorem 2.4, which implies
$$\mu(\{\omega\in \Omega ; |Z(n, x)(\omega)-\alpha_{n, x}|\ge \delta\})$$
$$=\mu\left (\left \{\omega\in \Omega ; \left |Z(n, x)(\omega)-(C)\int_{\Omega}Z(n, x)d \mu\right |\ge \delta\right \}\right )$$
$$\le \frac{(C)\int_{\Omega}(Z(n, x)-(C)\int_{\Omega}Z(n, x)d \mu)^{2}d \mu}{\delta^{2}}
=\frac{VAR_{ch}(Z(n, x))}{\delta^{2}}=\frac{\sigma_{n, x}^{2}}{\delta^{2}}.$$
$\hfill \square$

{\bf Remark 3.2.} Analyzing the proof of Theorem 3.1, it easily follows that without any change in its proof, the construction of the operators
$L_{n}(f)(x)$ can be slightly generalized by considering that not just $Z$ depends on $n$ and $x$, but also that $\mu$ may depend on $n$ and $x$ and $\Omega$ may depend on $n$.
More exactly, we can consider $L_{n}(f)(x)$ of the more general form
$$L_{n}(f)(x):=(C)\int_{\mathbb{R}}f(t) d\mu_{Z(n, x)}(t)=(C)\int_{\Omega_{n}}f \circ Z(n, x) d\mu_{n, x}, \, x\in I,\, n\in \mathbb{N},$$
where $\mu_{n, x}:{\cal{C}}_{n}\to [0, +\infty)$, $(n, x)\in \mathbb{N}\times I$, is a family of normalized, monotone and subadditive set functions on $\Omega_{n}$ and $\mu_{Z(n, x)}$ denotes the distribution of $Z(n, x)$ with respect to $\mu_{n, x}$.
This remark is useful in constructing several concrete examples of such operators.

Also, it is worth noting here that if we suppose that $\mu(\{\omega\in \Omega; Z(n, x)(\omega)\in I\}=1$, then the operators $L_{n}$ in the statement of Theorem 3.1 can be attached to continuous and bounded functions defined on a strict subinterval $I\subset \mathbb{R}$, $f:I\to \mathbb{R}$.

Indeed, let us extend $f$ to a function continuous and bounded, $f^{*}:\mathbb{R}\to \mathbb{R}$. Since $\mu$ is normalized, monotone and subadditive, this implies that $\mu_{Z(n, x)}$ is normalized, monotone and subadditive (see Definition 2.1, (iii) and the proof of Theorem 3.1).
But, for any monotone and subaddtive set function $\mu$ and $A$, $B$ with $\mu(B)=0$, we have
$\mu(A)\le \mu(A\bigcup B)\le \mu(A)+\mu(B)=\mu(A)$, which for any bounded function $F:A\bigcup B\to \mathbb{R}$ implies
$$(C)\int_{A\bigcup B}F d \mu=\int_{0}^{\infty}\mu(\{t\in A\bigcup B ; F(t)\ge \alpha\})d \alpha$$
$$-\int_{-\infty}^{0}[\mu(\{t\in A\bigcup B ; F(t)\ge \alpha\}) - \mu(A\bigcup B)]d \alpha$$
$$=\int_{0}^{\infty}\mu(\{t\in A ; F(t)\ge \alpha\})d \alpha-\int_{-\infty}^{0}[\mu(\{t\in A ; F(t)\ge \alpha\}) - \mu(A)]d \alpha$$
$$=(C)\int_{A}F d\mu.$$
Applying this for $A=I$, $B=\mathbb{R}\setminus I$, $F=f^{*}$ and $\mu=\mu_{Z(n, x)}$, we immediately get
$$(C)\int_{\mathbb{R}}f^{*}d \mu_{Z(n, x)}=(C)\int_{I}f d\mu_{Z(n, x)}.$$
In the next considerations, since as it is stated in Introduction, formula (\ref{int2}) does not hold in general, different approximation operators can directly be defined by the right-hand side in (\ref{int2}). In this sense, we present the following quantitative results.

{\bf Theorem 3.3.} {\it Denoting by ${\cal{P}}(\mathbb{R})$ the class of all subsets of $\mathbb{R}$, let $(\mathbb{R}, {\cal{C}})$ be a measurable space with ${\cal{C}}\subset {\cal{P}}(\mathbb{R})$ and
$\mu:{\cal{C}}\to [0, +\infty)$, be a monotone and submodular set function.

For $\lambda_{n, x}:\mathbb{R}\to \mathbb{R}_{+}$, $n\in \mathbb{N}$, $x\in \mathbb{R}$, Choquet densities with respect to $\mu$, (that is, $(C)\int_{\mathbb{R}}\lambda_{n, x}(t) d \mu(t) = 1$), let us define by $UC(\mathbb{R})$, the class of all functions $f:\mathbb{R}\to \mathbb{R}_{+}$, uniformly continuous on $\mathbb{R}$, such that $f\cdot \lambda_{n, x}$ are ${\cal{C}}$-measurable and $T_{n}(f)(x)<+\infty$, for all $n\in \mathbb{N}$, $x\in \mathbb{R}$, where
$$T_{n}(f)(x)=(C)\int_{\mathbb{R}}f(t)\cdot \lambda_{n, x}(t) d \mu(t).$$

Then, denoting $\varphi_{x}(t)=|t-x|$,
for all $x\in \mathbb{R}$, $n\in \mathbb{N}$ and $\delta >0$ we have
$$|T_{n}(f)(x)-f(x)|\le \left [1+\frac{T_{n}(\varphi_{x})(x)}{\delta}\right ]\cdot \omega_{1}(f ; \delta)_{\mathbb{R}}.$$}

{\bf Proof.} According to Remark 2.3, (i), we immediately get that each $T_{n}$ is a monotone, subadditive, positive homogenous operator on the space $UC(\mathbb{R})$. It is also worth noting that if $f, g\in UC(\mathbb{R})$ then $f+g\in UC(\mathbb{R})$ and $\alpha \cdot f\in UC(\mathbb{R})$, for $\alpha\ge 0$.

Let $f,g\in UC(\mathbb{R})$. We have $f=f-g+g\leq |f-g|+g$, which successively implies $T_{n}(f)(x)\leq
T_{n}(|f-g|)(x)+T_{n}(g)(x)$, that is $T_{n}(f)(x)-T_{n}(g)(x)\leq
T_{n}(|f-g|)(x)$.

Writing now $g=g-f+f\leq|f-g|+f$ and applying the above reasonings, it
follows $T_{n}(g)(x)-T_{n}(f)(x)\leq T_{n}(|f-g|)(x)$, which combined with
the above inequality gives $|T_{n}(f)(x)-T_{n}(g)(x)|\leq T_{n}(|f-g|)(x)$.

Then, from the identity
\[
T_{n}(f)(x)-f(x)=[T_{n}(f)(x)-f(x)\cdot
T_{n}(e_{0})(x)]+f(x)[T_{n}(e_{0})(x)-1],
\]
by using the above inequality too, it follows
\[
|f(x)-T_{n}(f)(x)|\leq|T_{n}(f(x))(x)-T_{n}(f(t))(x)|+|f(x)|%
\cdot|T_{n}(e_{0})(x)-1|
\]
\[
\leq T_{n}(|f(t)-f(x)|)(x)+|f(x)|\cdot|T_{n}(e_{0})(x)-1|.
\]
Now, since for all $t,x\in I$ we have
\[
|f(t)-f(x)|\leq\omega_{1}(f;|t-x|)_{I}\leq\left[ \frac{1}{\delta }|t-x|+1%
\right] \omega_{1}(f;\delta)_{I},
\]
replacing above and taking into account that $T_{n}(e_{0})(x)=e_{0}(x)$ (here $e_{0}(x)=1$, for all $x\in \mathbb{R}$), we immediately obtain that for all $n\in \mathbb{N}$, $x\in \mathbb{R}$ and $\delta >0$, we have
$$|T_{n}(f)(x)-f(x)|\le \left [1+\frac{1}{\delta}T_{n}(\varphi_{x})(x)\right ]\omega_{1}(f; \delta)_{\mathbb{R}}.$$
$\hfill \square$

{\bf Remark 3.4.} An important problem in Theorem 3.3 is to determine the functions $f$ with $T_{n}(f)(x)<+\infty$, for all $n\in \mathbb{N}$, $x\in \mathbb{R}$. Since $T_{n}(e_{0})(x)=1$, from the positive homogeneity of $T_{n}$ it easily follows that for $f_{0}(t)=c>0$ for all $t\in \mathbb{R}$, we have $T_{n}(f_{0})(x)=c$, for all $x\in \mathbb{R}$, $n\in \mathbb{R}$. Then, for any bounded $f:\mathbb{R}\to \mathbb{R}_{+}$ and $\alpha \ge 0$, denoting $\|f\|=\sup\{f(x) ; x\in \mathbb{R}\}<+\infty$, we get $\{t\in \mathbb{R} ; f(t)\cdot \lambda_{n, x}(t)\ge \alpha \}\subset \{t\in \mathbb{R} ; \|f\|\cdot \lambda_{n, x}(t)\ge \alpha\}$, which implies
$$\mu(\{t\in \mathbb{R} ; f(t)\cdot \lambda_{n, x}(t)\ge \alpha \})\le \mu(\{t\in \mathbb{R} ; \|f\|\cdot \lambda_{n, x}(t)\ge \alpha\})$$
and therefore
$$T_{n}(f)(x)=\int_{0}^{+\infty}\mu(\{t\in \mathbb{R} ; f(t)\cdot \lambda_{n, x}(t)\ge \alpha \})d \alpha$$
$$\le \int_{0}^{+\infty}\mu(\{t\in \mathbb{R} ; \|f\|\cdot \lambda_{n, x}(t)\ge \alpha \})d \alpha=\|f\|\cdot T_{n}(e_{0})(x)=\|f\|<+\infty.$$
Also, from these reasonings, it follows that if for an unbounded function $F_{0}:\mathbb{R}\to \mathbb{R}_{+}$ we have $T_{n}(F_{0})(x)<+\infty$, for all $x\in \mathbb{R}$, $n\in \mathbb{N}$, then for any unbounded function $f$ satisfying $f(t)\le F_{0}(t)$, for all $t\in \mathbb{R}$, we
have $T_{n}(f)(x)<+\infty$, for all $x\in \mathbb{R}$, $n\in \mathbb{N}$.

{\bf Remark 3.5.} Analyzing the proof of Theorem 3.3, it is clear that $\mu$ may depend on $n$ and $x$ too. Also, it is worth noting that the operators studied by Theorem 3.1 (and those more general defined in Remark 3.2), in general do not coincide with those
studied by Theorem 3.3, fact which implies that these two theorems represent distinct results. This is due to the fact that for a monotone and subadditive Choquet set function $\mu$, in general, we don't have the equality
$$(C)\int_{\mathbb{R}}f(t) d\mu_{Z(n, x)}(t)=(C)\int_{\mathbb{R}}f(t) \lambda_{n, x}(t)d \mu(t), \mbox{ for all } f,$$
where $\lambda_{n, x}:\mathbb{R}\to \mathbb{R}_{+}$ represents the Choquet density of $Z(n, x)$, as defined by Definition 2.1, (iii).

However, as it was pointed out in Introduction too, due to a Radon-Nikodym result for the Choquet integral in \cite{Graf} (see also \cite{Mol}, p. 75, Theorem 5.10), for very special subclasses of monotone and subadditive set functions, the formula for $L_{n}(f)(x)$ in (\ref{int2}) still holds.  More exactly, if $\mu$ and $\nu$ are monotone, subadditive and continuous from below set functions, where $\nu(A)=(C)\int_{A}f d\mu$, for all $A\in {\cal{B}}$, if the couple $(\mu, \nu)$ has a strong decomposition property and if $\mu(A)=0$ implies $\nu(A)=0$, then the formula (\ref{int2}) holds for all non negative functions $f$.

\section{Concrete Examples}

In this section, the results in Section 3 are illustrated by several examples.

{\bf Example 4.1.} Here we consider approximation operators generated by the construction in Theorem 3.1 (see also Remark 3.2 too).

Take $\Omega_{n}=\{0, 1, 2, ..., n\}$ and consider $Z(n, x):\Omega_{n}\to [0, 1]$ defined by $Z(n, x)(k)=\frac{k}{n}$, for all $k\in \Omega_{n}$ (in fact, here $Z$ depends only on $n$).

Firstly, if we define $\mu_{n, x}:{\cal{P}}(\Omega_{n})\to \mathbb{R}_{+}$ by
$\mu_{n, x}(\{i\})=p_{n, i}(x)={n\choose i}x^{i}(1-x)^{n-i}$, $i\in \{0, 1, ..., n\}$, $x\in [0, 1]$ and for any $A\subset \Omega_{n}$,
$\mu_{n, x}(A)=\sum_{i\in A}p_{n, i}(x)$, then $\mu_{n, x}$ is normalized, monotone and additive on ${\cal{P}}(\Omega_{n})$. Also, defining  $\mu_{Z(n, x)}$ as in Definition 2.1, (iii), by the last formula of calculation in Definition 2.1, (ii), for $L_{n}(f)(x)$ in Theorem 3.1 (taking into account Remark 3.2 too), we easily recapture the classical Bernstein polynomials
$$L_{n}(f)(x)=(C)\int_{\Omega_{n}}f\circ Z(n, x) d\mu_{n, x}=\sum_{i=0}^{n}f(i/n)p_{n, i}(x)=B_{n}(f)(x),\, x\in [0, 1].$$

Secondly, define $\mu_{n, x}:\Omega_{n}\to [0, 1]$, $n\in \mathbb{N}$, $n\ge 2$, $x\in [0, 1]$, by
$\mu_{n, x}(\{i\})=p_{n, i}(x)$, if $i\in \Omega_{n}\setminus\{1\}$, $\mu(\{1\})=\varphi_{n, 1}(x)$, where $\varphi_{n, 1}(x)$ is chosen arbitrary  satisfying
$p_{n, 1}(x)\le \varphi_{n, 1}(x)\le p_{n, 1}(x)+\min\{p_{n, i}(x) ; i\in \Omega_{n}\setminus \{1\}\}$, for all $x\in [0, 1]$, $n\in \mathbb{N}$, $n\ge 2$ and
for any $A\subset \Omega_{n}$, $A\not = \Omega_{n}$, $\mu_{n, x}(A)=\sum_{i\in A}p_{n, i}(x)$ if $1\not\in A$, $\mu_{n, x}(A)=\sum_{i\in A, i\not=1}p_{n, i}(x)+\varphi_{n, 1}(x)$, if $1\in A$. Finally, define $\mu_{n, x}(\Omega_{n})=1$

It easily follows that $\mu_{n, x}$ is a normalized, monotone and non-additive but subadditive set function.

Now, let us consider that $f:[0, 1]\to \mathbb{R}$ is non-decreasing on $[0, 1]$ and  $Z(n, x)(i)=\frac{i}{n}$, $i\in \Omega_{n}$. By the last formula of calculation in Definition 2.1, (ii), since $f \circ Z(n, x)$ is non-decreasing, we immediately get the Bernstein-Choquet kind operator
$$L_{n}(f)(x)=f(0)[\mu_{n, x}(\{0, 1, ..., n\})-\mu_{n, x}(\{1, 2, ..., n\})]$$
$$+f(1/n)[\mu_{n, x}(\{1, 2, .., n\})-\mu_{n, x}(\{2, 3, ..., n\})]$$
$$+\sum_{i=2}^{n}f(i/n)[\mu_{n, x}(\{i, ..., n\})-\mu_{n, x}(\{i+1, ..., n\})]$$
$$=f(0)[p_{n, 0}(x)+p_{n, 1}(x)-\varphi_{n, 1}(x)]+f(1/n)\varphi_{n, 1}(x)+
\sum_{i=2}^{n}f(i/n)p_{n, i}(x)$$
$$=B_{n}(f)(x)+f(0)\cdot [p_{n, 1}(x)-\varphi_{n, 1}(x)]+f(1/n)\cdot [\varphi_{n, 1}(x)-p_{n, 1}(x)],$$
that is
\begin{equation}\label{int3}
L_{n}(f)(x)=B_{n}(f)(x)+[f(1/n)-f(0)]\cdot [\varphi_{n, 1}(x)-p_{n, 1}(x)],
\end{equation}
where $\varphi_{n, 1}(x)$ satisfies
\begin{equation}\label{int4}
0\le \varphi_{n, 1}(x)-p_{n, 1}(x)\le \min\{p_{n, 0}(x), p_{n, n}(x)\}=\min\{x^{n}, (1-x)^{n}\}\le \left (\frac{1}{2}\right )^{n},
\end{equation}
which implies that $\varphi_{n, 1}-p_{n, 1}\to 0$, uniformly on $[0, 1]$.

This immediately implies that
$$L_{n}(e_{1})(x)=B_{n}(e_{1})(x)+[f(1/n)-f(0)]\cdot [\varphi_{n, 1}(x)-p_{n, 1}(x)]\to x,$$ (as $n\to \infty$), uniformly on
$[0, 1]$.

Similarly, if $f$ is non-increasing on $[0, 1]$, since $f\circ Z(n, x)$ is non-increasing on $[0, 1]$, then by the formula of calculation in Definition 2.1, (ii), we easily get that
$$L_{n}(f)(x)=(C)\int_{\Omega_{n}}f\circ Z(n, x) d \mu_{n, x}$$
$$=B_{n}(f)(x)+[f(1/n)-f(1)][\varphi_{n, 1}(x)-p_{n, n}(x)],$$
formula which is somehow symmetric to the formula in the case when $f$ is non-decreasing on $[0, 1]$.

Now, in order to evaluate the order of the expression $\sigma_{n, x}^{2}$ in the statement of Theorem 3.1, let us denote
$$Y(n, x)(t)=(Z(n, x)(t)-L_{n}(e_{1})(x))^{2}=(x-Z(n, x)(t)+\frac{1}{n}[\varphi_{n, 1}(x)-p_{n, 1}(x)])^{2}.$$
We clearly have
$$Y(n, x)(t)\le 2\cdot  \left [(x-Z(n, x)(t))^{2}+\frac{1}{n^{2}}\cdot (\varphi_{n, 1}(x)-p_{n, 1}(x))^{2}\right ],$$
which implies
$$\sigma_{n, x}^{2}=L_{n}(Y(n, x)(\cdot))(x)$$
$$\le 2 \cdot (C)\int_{\Omega_{n}}\left [(x-Z(n, x)(t))^{2}d \mu_{n, x}(t)+\frac{1}{n^{2}}(\varphi_{n, 1}(x)-p_{n, 1}(x))^{2}\right ] d \mu_{n, x}(t)$$
$$\le 2 \cdot (C)\int_{\Omega_{n}}\left [(x-Z(n, x)(t))^{2}+\frac{1}{n^{2}} \right ]d\mu_{n, x}(t)$$
$$=2 \cdot (C)\int_{\Omega_{n}}(x-Z(n, x)(t))^{2}d\mu_{n, x}(t)+\frac{2}{n^{2}}.$$
Denoting $W(n, x)(t)=(x-Z(n, x)(t))^{2}$, it remains to put in increasing order the values $W(n, x)(k/n)$, $k=0, 1, ..., n$, depending of course on the values of $x$.

Let $k_{0}\in \{0, ..., n-1\}$ be fixed. For $x\in \left [\frac{k_{0}}{n}, \frac{k_{0}+1}{n}\right ]$, by simple calculation we get :
if $x\in \left [\frac{k_{0}}{n}, \frac{2k_{0}+1}{2 n}\right ]$, then
$$W(n, x)(k_{0}/n)\le W(n, x)((k_{0}+1)/n)\le W(n, x)((k_{0}-1)/n)$$
$$\le W(n, x)((k_{0}+2)/n)\le W(n, x)((k_{0}-2)/n)\le W(n, x)((k_{0}+3)/n)$$
$$\le ...\le W(n, x)(2k_{0}/n)\le W(n, x)(0/n)\le W(n, x)((2k_{0}+1)/n)$$
$$\le W(n, x)((2k_{0}+2)/n)\le ... \le W(n, x)(1),$$
and if $x\in \left [\frac{2k_{0}+1}{2 n}, \frac{k_{0}+1}{n}\right ]$, then
$$W(n, x)((k_{0}+1)/n)\le W(n, x)(k_{0}/n)\le W(n, x)((k_{0}+2)/n)$$
$$\le W(n, x)((k_{0}-1)/n)\le W(n, x)((k_{0}+3)/n)\le W(n, x)((k_{0}-2)/n)$$
$$\le ...\le W(n, x)((2k_{0}+1)/n)\le W(n, x)(0/n)\le W(n, x)((2k_{0}+2)/n)$$
$$\le W(n, x)((2 k_{0}+3)/n)\le ... \le W(n, x)(1).$$
For example, if we take $k_{0}=0$ and $x\in [0, 1/(2n)]$, then we easily get
$$W(n, x)(0)\le W(n, x)(1/n)\le W(n, x)(2/n)\le ... \le W(n, x)(1),$$
and therefore, by the formula in Definition 2.1, (ii), we easily obtain
$$L_{n}(W(n, x)(\cdot))(x)=x^{2}[p_{n, 0}(x)+p_{n, 1}(x)-\varphi_{n, 1}(x)]+(x-1/n)^{2}\cdot \varphi_{n, 1}(x)$$
$$+\sum_{i=2}^{n}(x-i/n)^{2}p_{n, i}(x)=B_{n}((x-\cdot )^{2})(x)+[\varphi_{n, 1}(x)-p_{n, 1}(x)]\cdot [(x-1/n)^{2}-x^{2}]$$
$$=\frac{x(1-x)}{n}+[\varphi_{n, 1}(x)-p_{n, 1}(x)]\cdot [1/n^{2} - 2x/n],$$
where $0\le \varphi_{n, 1}-p_{n, 1}(x)\le \frac{1}{2^{n}}$, for all $x\in [0, 1]$, $n\ge 2$.

For another example, consider $k_{0}=2$ and $x\in [2/n, 5/(2n)]$. Then, by formula in Definition 2.1, (ii), we can write
$$L_{n}(W(n, x)(\cdot))(x)=(x-2/n)^{2}[p_{n, 0}(x)+p_{n, 1}(x)-\varphi_{n, 1}(x)]+(x-3/n)^{2}\cdot \varphi_{n, 1}(x)$$
$$+(x-1/n)^{2}\cdot p_{n, 2}(x)+(x-4/n)^{2}\cdot p_{n, 3}(x)+x^{2}\cdot p_{n, 4}(x)+\sum_{i=5}^{n}(x-i/n)^{2}\cdot p_{n, i}(x)$$
$$=B_{n}((x-\cdot )^{2})(x)+x^{2}[p_{n, 4}(x)-p_{n, 0}(x)]+(x-1/n)^{2}[p_{n, 3}(x)-p_{n, 1}(x)]$$
$$+(x-2/n)^{2}[p_{n, 0}(x)+p_{n, 1}(x)-\varphi_{n, 1}(x)-p_{n, 2}(x)]+(x-3/n)^{2}[\varphi_{n, 1}(x)-p_{n, 3}(x)]$$
$$+(x-4/n)^{2}[p_{n, 3}(x)-p_{n, 4}(x)],$$
which for $x\in [2/n, 5/(2n)]$ immediately implies
$$|L_{n}(W(n, x)(\cdot))(x)|\le \frac{C}{n^{2}},$$
where $C>0$ satisfies an inequality of the form $C\le C_{0}\cdot n$ on each compact subinterval of $[0, 1]$, with $C_{0}$ an absolute constant.

Summarizing all the cases, we easily get that $L_{n}(W(n, x)(\cdot))(x)\to 0$, uniformly on each compact subinterval of $[0, 1]$, which by Theorem 3.1 implies that
$L_{n}(f)\to f$, uniformly on each compact subinterval of $[0, 1]$, for any continuous function $f:[0, 1]\to \mathbb{R}$.

It is worth noting that for continuous non-decreasing functions $f$, the order of approximation given by $L_{n}(f)(x)$ in (\ref{int3}) is not worst than the order of approximation given by the Bernstein polynomials $B_{n}(f)(x)$. However, for example, in the subclass of increasing and concave functions on $[0, 1]$, $L_{n}(f)(x)$ approximates better than the Bernstein polynomials $B_{n}(f)(x)$.
Indeed, since for all $x\in (0, 1)$ and $n\in \mathbb{N}$ we have $B_{n}(f)(x)-f(x)<0$ and $[f(1/n)-f(0)]\cdot [\varphi_{n, 1}(x)-p_{n, 1}(x)]>0$,
by (\ref{int3}) we easily obtain
$$|L_{n}(f)(x)-f(x)| < \max\{|B_{n}(f)(x)-f(x)|, [f(1/n)-f(0)]\cdot [\varphi_{n, 1}(x)-p_{n, 1}(x)]\},$$
for all $x\in (0, 1)$ and $n\in \mathbb{N}$, which by (\ref{int4}) easily leads us to the above conclusion.

A similar approach can be made for the subclass of decreasing and concave functions on $[0, 1]$.

Other $n$ examples of nonlinear Bernstein-Choquet type operators on $C[0, 1]$ which satisfy the approximation properties in Theorem 3.1, can be
obtained if in the definition of $\mu_{n, x}$ we replace the index $1\in \{0, ..., n\}$ by another arbitrary fixed $i_{0}\in \{0, ..., n\}$ (that is
$\mu_{n, x}(\{i_{0}\})=\varphi_{n, i_{0}}(x)$, $\mu_{n, x}(\{i\})=p_{n, i}(x)$ if $i\not=i_{0}$) and repeat the same construction as above.
Also, obviously that these kinds of constructions can be applied to other kinds of Bernstein-type operators too.

{\bf Example 4.2.} In what follows, we consider an example of operator satisfying Theorem 3.3.

Let us define the nonlinear Picard-Choquet operators attached to $f:\mathbb{R}\to \mathbb{R}_{+}$ and to the monotone, submodular set function
$\mu$, by
$$T_{n}(f)(x)=\frac{1}{c(n, x)} \cdot (C)\int_{\mathbb{R}}f(t)e^{-n|x-t|}d \mu(t)$$
$$=\frac{1}{c(n, x)} \cdot \int_{0}^{\infty}\mu[F_{\alpha}(f(\cdot) e^{-n|x-\cdot|})]d\alpha,$$
where $c(n, x)=\int_{0}^{\infty}\mu[F_{\alpha}(e^{-n|x-\cdot|})]d \alpha$, $F_{\alpha}(e^{-n|x-\cdot|})=\{t\in \mathbb{R} : e^{-n|x-t|}\ge \alpha\}$.
By simple calculation we get $F_{\alpha}(e^{-n|x-\cdot|})=\emptyset$ for $\alpha >1$ and $F_{\alpha}(e^{-n|x-\cdot|})=\left [x+\frac{ln(\alpha)}{n}, x-\frac{ln(\alpha)}{n}\right ]$, if $\alpha\le 1$, which leads us to
$$c(n, x)=\int_{0}^{1}\mu\left (\left [x+\frac{ln(\alpha)}{n}, x-\frac{ln(\alpha)}{n}\right ]\right )d\alpha.$$
Indeed, $\{t\in \mathbb{R} ; e^{-n|t-x|}\ge \alpha\}=\emptyset$ for $\alpha>1$ and for all $0\le \alpha\le 1$ we have
$$\{t\in \mathbb{R} ; e^{-n|t-x|}\ge \alpha\}$$
$$=\{t\in \mathbb{R} ; t\ge x, e^{-n(t-x)}\ge \alpha\}\bigcup \{t\in \mathbb{R} ; t < x, e^{-n(x-t)}\ge \alpha\}$$
$$\{t\in \mathbb{R} ;x\le t\le \frac{nx-ln(\alpha)}{n}\}\bigcup \{t\in \mathbb{R}; \frac{ln(\alpha)+n x}{n}\le t\le x\}$$
$$=\left [\frac{nx+ln(\alpha)}{n}, \frac{nx-ln(\alpha)}{n}\right ].$$

By Theorem 3.3, for all $n\in \mathbb{N}$, $x\in \mathbb{R}$ and $\delta >0$ we get
\begin{equation}\label{Pici}
|T_{n}(f)(x)-f(x)|\le \left [1+\frac{1}{\delta}T_{n}(\varphi_{x})(x)\right ]\omega_{1}(f; \delta)_{\mathbb{R}},
\end{equation}
where $\varphi_{x}(t)=|t-x|$. Therefore, the convergence of $T_{n}(f)$ to $f$ one relies on the convergence to zero, as $n\to \infty$, of the
quantity
$$T_{n}(\varphi_{x})(x)=\frac{1}{c(n, x)}\cdot (C)\int_{\mathbb{R}}|x-t|e^{-n|x-t|}d\mu(t)$$
$$=\frac{1}{c(n, x)}\cdot \int_{0}^{\infty}\mu[\{t\in \mathbb{R} ; |t-x|\cdot e^{-n|t-x|}\ge \alpha\}]d \alpha.$$

Since any possibility measure also is monotone and submodular (see Remark 2.3, (ii)), let use consider the possibilistic measure
$\mu_{n, x}(A)=\sup\{e^{-n|t-x|} ; t\in A\}, \mbox{ if } A\subset \mathbb{R}, A\not=\emptyset$ and
$\mu_{n, x}(\emptyset)=\sup\{e^{-n|t-x|} ; t\in \emptyset\}=0$.
The Picard-Choquet integral for $f:\mathbb{R}\to \mathbb{R}_{+}$ becomes
$$T_{n}(f)(x)=\frac{1}{c(n, x)}\cdot (C)\int_{\mathbb{R}}f(t)\cdot e^{-n|t-x|} d\mu_{n, x}(t)$$
$$=\frac{1}{c(n, x)}\cdot \int_{0}^{\infty}\sup\{e^{-n|t-x|} ; t\in \mathbb{R}, f(t)\cdot e^{-n|t-x|}\ge \alpha\}d \alpha$$
and its convergence to $f$ depends on the upper estimates of the quantity
$$T_{n}(\varphi_{x})(x)=\frac{1}{c(n, x)}\cdot \int_{0}^{\infty}\sup\{e^{-n|t-x|} ; t\in \mathbb{R}, |t-x|e^{-n|t-x|}\ge \alpha\}d \alpha.$$
By the above formula of $c(n, x)$ with $\mu$ replaced by the possibilistic measure $\mu_{n, x}$, we obtain
$$c(n, x)=\int_{0}^{\infty}\sup\{e^{-n|t-x|} ; t\in \mathbb{R}, e^{-n|t-x|}\ge \alpha \}d \alpha$$
$$=\int_{0}^{1}\sup\{e^{-n|t-x|} ; t\in [x+\ln(\alpha)/n, x-\ln(\alpha)/n]\}d\alpha.$$
Since $t\in [x+\ln(\alpha)/n, x-\ln(\alpha)/n]$ is equivalent with $0\le |t-x|\le -\frac{\ln(\alpha)}{n}$, it follows
$$c(n, x)=\int_{0}^{1}\sup\{e^{-n|t-x|} ; t\in [x+\ln(\alpha)/n, x-\ln(\alpha)/n]\}d\alpha=\int_{0}^{1}1 d\alpha=1.$$
On the other hand, we have
$$T_{n}(\varphi_{x})(x)=\frac{1}{c(n, x)}\cdot \int_{0}^{\infty}\sup\{e^{-n|t-x|} ; t\in \mathbb{R}, |t-x|e^{-n|t-x|}\ge \alpha\}d \alpha$$
$$=\int_{0}^{\infty}\sup\{e^{-ny} ; y\ge 0, y e^{-ny}\ge \alpha\}d \alpha.$$
Now, denoting $F(v)=v e^{-n v}$, $v\ge 0$, we have $F^{\prime}(v)=(1-n v)e^{-n v}$, which immediately implies that $v=\frac{1}{n}$ is a maximum point for $F$ on $[0, +\infty)$ and $F(1/n)=\frac{1}{n e}$ is the maximum value for $F$.

This implies that for $\alpha > 1/(n e)$ we have $\{t\in \mathbb{R} ; |t-x|\cdot e^{-n|t-x|}\ge \alpha\}=\emptyset$ and therefore it follows
$$T_{n}(\varphi_{x})(x)=\int_{0}^{1/(n e)}\sup\{e^{-ny} ; y\ge 0, y e^{-ny}\ge \alpha\}d \alpha\le \int_{0}^{1/(n e)}1\cdot d\alpha=\frac{1}{n e}.$$
Then, choosing $\delta=\frac{1}{n e}$ in Theorem 3.3, we immediately get the good approximation estimate
$$|T_{n}(f)(x)-f(x)|\le 2\omega_{1}\left (f; \frac{1}{n e}\right )_{\mathbb{R}}.$$
It is worth mentioning that however, there are classes of functions for which $T_{n}(f)(x)$ gives an essentially better estimate. For example,
for all functions of the form $f(t)=ce^{-\lambda t}$, with $c$ a real constant and $\lambda >0$, we have $T_{n}(f)(x)=f(x)$,
for all $x\in \mathbb{R}$ and $n > \lambda$. For the simplicity of calculation, take, for example, $f(t)=e^{-t}$.  We get
$$T_{n}(f)(x)=\int_{0}^{+\infty}\sup\{e^{-n|t-x|} ; t\in \mathbb{R}, e^{-t}\cdot e^{-n|t-x|}\ge \alpha\}d \alpha$$
$$=\int_{0}^{1}\sup\{e^{-n|t-x|} ; t\in \mathbb{R}, e^{-t}\cdot e^{-n|t-x|}\ge \alpha\}d \alpha.$$
But for $\alpha \in [0, 1]$ we have
$$\{t\in \mathbb{R}, e^{-t}\cdot e^{-n|t-x|}\ge \alpha\}$$
$$=\{t\in \mathbb{R}, t\ge x, e^{-t}\cdot e^{-n|t-x|}\ge \alpha\}\bigcup \{t\in \mathbb{R}, t<x, e^{-t}\cdot e^{-n|t-x|}\ge \alpha\}$$
$$=\{t\in \mathbb{R}, x\le t \le \frac{n x -\ln(\alpha)}{n+1}\}\bigcup \{t\in \mathbb{R}, \frac{n x + \ln(\alpha)}{n-1}\le t < x\}.$$
But $x\le \frac{n x-\ln(\alpha)}{n+1}$ if and only if $\alpha \le e^{-x}$, and $\frac{n x + \ln(\alpha)}{n-1}\le  x$ if and only if $\alpha \le e^{-x}$, which immediately implies
$$T_{n}(f)(x)=\int_{0}^{e^{-x}}\sup\left \{e^{-n|t-x|} ; t\in \left [\frac{n x + \ln(\alpha)}{n-1}, \frac{n x-\ln(\alpha)}{n+1}\right ]\right \}d\alpha$$
$$=\int_{0}^{e^{-x}}1\cdot d\alpha=e^{-x},$$
since $x\in [\frac{n x + \ln(\alpha)}{n-1}, \frac{n x-\ln(\alpha)}{n+1}]$.

By simple calculation it is easy to check that for the classical Picard operator $P_{n}(f)(x)=\frac{n}{2}\cdot \int_{\mathbb{R}}f(t)\cdot e^{-n|t-x|}d t$, we don't have $P_{n}(f)(x)=f(x), x\in \mathbb{R}$, for $f(x)=e^{-x}$.

Now, instead of the measure of possibility in the definition of the Picard-Choquet operator we can consider $\mu:{\cal{M}}(\mathbb{R})\to \mathbb{R}_{+}$
given by $\mu(A)=\sqrt{m(A)}$, where ${\cal{M}}(\mathbb{R})$ denotes the class of all bounded, Lebesgue measurable subsets of $\mathbb{R}$
and $m(A)$ denotes the Lebesgue measure. Since $\gamma(t)=\sqrt{t}$, $t\ge0$, is increasing and concave, according to Remark 2.3, (ii), it follows that $\mu$ is a monotone and submodular set function on ${\cal{M}}(\mathbb{R})$.

In this case, we easily get
$$c(n, x)=\frac{\sqrt{2}}{\sqrt{n}}\cdot \int_{0}^{1}\sqrt{-\ln(\alpha)} d\alpha =\frac{\sqrt{2}}{\sqrt{n}}\cdot \int^{+\infty}_{0}t^{1/2}e^{-t} dt
=\frac{\sqrt{2}}{\sqrt{n}}\cdot \Gamma(3/2)=\frac{\sqrt{\pi}}{\sqrt{2} \sqrt{n}}$$
and by similar reasonings with those for the possibility measure, we can find an estimate for $T_{n}(\varphi_{x})(x)$ and consequently, by Theorem 3.3 a quantitative estimate in approximation of $f$ by $T_{n}(f)$.

As it was pointed out by Remark 3.4 too, the above estimates hold for any bounded and uniformly continuous function $f:\mathbb{R}\to \mathbb{R}_{+}$.
But it can easily be determined unbounded, uniformly continuous functions $F_{0}$, such that $T_{n}(F_{0})(x)<+\infty$, for all $x\in \mathbb{R}$, $n\in \mathbb{N}$, and for which the above estimate still holds. Then, according to Remark 3.4 we get that for all $f$ satisfying $0\le f(t)\le F_{0}(t)$, for all $t\in \mathbb{R}$, we have $T_{n}(f)(x)<+\infty$, $n\in \mathbb{N}$, $x\in \mathbb{R}$.

{\bf Example 4.3.} Similar to the Picard-Choquet operators, the Gauss-We\-ier\-strass-Choquet operators given by
$$W_{n}(f)(x)=\frac{1}{c(n, x)} \cdot (C)\int_{\mathbb{R}}f(\xi)e^{-n(x-\xi)^{2}}d \mu(\xi)$$
could be studied, with $c(n, x)=(C)\int_{\mathbb{R}}e^{-n(x-\xi)^{2}}d \mu(\xi)$,
$n\in \mathbb{N}$, $f:\mathbb{R}\to \mathbb{R}_{+}$. In this case,  we get
$$c(n, x)=\int_{0}^{1}\mu\left (\left [x-\sqrt{-ln(\alpha)/n}, x+\sqrt{-ln(\alpha)/n}\right ]\right )d\alpha,$$
and the convergence of $W_{n}(f)$ to $f$ (as $n\to 0$) one relies on the convergence to zero of the
quantity $W_{n}(\varphi_{x})(x)=\frac{1}{c(n, x)} \cdot (C)\int_{\mathbb{R}}|x-\xi|e^{-n(x-\xi)^{2}}d \mu(\xi)$.
As in the previous Example 4.2, for the possibility measure $\mu_{n, x}(A)=\sup\{e^{-n|t-x|} ; t\in A\}, \mbox{ if } A\subset \mathbb{R}, A\not=\emptyset$ and for $\mu(A)=\sqrt{m(A)}$, where $"m"$ is the Lebesgue measure), we can derive that $W_{n}(\varphi_{x})(x)\to 0$ as $n\to \infty$ and obtain quantitative estimates in approximation by using Theorem 3.3.

\end{document}